\documentclass[11pt]{article}

\usepackage{amscd,amsmath,amsthm,latexsym}
\usepackage{amsfonts}
\usepackage{amssymb}
\usepackage{longtable}
\usepackage[all]{xy}
\usepackage{graphicx}
\usepackage{esint}
\usepackage{bbm}

\newenvironment{pf}{\medskip\noindent{\it Proof.}\enspace}
{\hfill \newline \smallskip}

\numberwithin{equation}{section}

\newcommand{\bP}{\mathbb{P}}

\newcommand{\bE}{{\mathbb{E}}}
\newcommand{\bR}{{\mathbb{R}}}

\newcommand{\af}{\alpha}
\newcommand{\bt}{\beta}
\newcommand{\gm}{\gamma}

\newcommand{\ep}{\varepsilon}
\newcommand{\zt}{\zeta}
\newcommand{\te}{\theta}
\newcommand{\ld}{\lambda}

\newcommand{\kp}{\kappa}

\newcommand{\wt}{\widetilde}
\newcommand{\sms}{\setminus}

\newcommand{\Gm}{\Gamma}
\newcommand{\Ld}{\Lambda}

\def\1{{\bf 1}}

 \def\sH {{\cal H}}

 \def\bE {{\mathbb E}}

 \def\bN {{\mathbb N}} 
\def\bP {{\mathbb P}}  \def\bR {{\mathbb R}}

\def\R {{\mathbb R}} \def\RR {{\mathbb R}}

\def\wt{\widetilde}
\def\pf{\noindent{\bf Proof.} }
\def\beq{\begin{equation}}
\def\eeq{\end{equation}}
\def\bee{\begin{equation}}
\def\eee{\end{equation}}
\def\ep{\epsilon}
\def\de{\delta}

\def\nn{\nonumber}

\newtheorem{thm}{Theorem}[section]
\newtheorem{lemma}[thm]{Lemma}

\newtheorem{definition}[thm]{Definition}

\newtheorem{cor}[thm]{Corollary}
\newtheorem{remark}[thm]{Remark}
\newtheorem{example}[thm]{Example}
\numberwithin{equation}{section}

\def\qed{{\hfill $\Box$ \bigskip}}

\def\LL{{\cal L}}

\def\RR{{\mathbb R}}
\def\R{{\mathbb R}}

\def\E{{\mathbb E}}

\def\P{{\mathbb P}}

\def\wt{\widetilde}
\def\pf{\noindent{\bf Proof.} }
\def\beq{\begin{equation}}
\def\eeq{\end{equation}}
\def\bee{\begin{equation}}
\def\eee{\end{equation}}
\def\rd{\mathrm d}
\def\Dc{{\overline D}^c}
\def\Uc{{\overline U}^c}
\def\Rd{\bR^d}
\def\RL{R_{\textrm{Lip}}}
\def\LL{\Lambda_{\textrm{Lip}}}

\begin{document}

\bibliographystyle{plain}

\title{\bf Tangential limits for harmonic functions with respect to $\phi(\Delta)$
 : stable and beyond
}

\author{{\bf Jaehoon Kang}\thanks{
This work was supported by the National Research Foundation of Korea(NRF) grant funded by the Korea government(MSIP) (No.2009-0083521)}
 \quad and  \quad {\bf Panki Kim}\thanks{This work was supported by the National Research Foundation of Korea(NRF) grant funded by the Korea government(MEST) (2013004822)
}}

\date{}

\maketitle

\begin{abstract}
In this paper, we discuss tangential limits for 
regular harmonic functions 
with respect to $\phi(\Delta):=-\phi(-\Delta)$
in the $C^{1,1}$ open set $D$ in $\R^d$, where $\phi$ is the complete Bernstein function and $d \ge 2$. 
When the exterior function $f$ is local $L^p$-H\"older continuous of order 
$\bt$ on $D^c$ with $ p\in(1,\infty]$ and $\bt>1/p$, 
for a large class of Bernstein function $\phi$,  
we show that the regular harmonic function $u_f$ with respect to $\phi(\Delta)$, 
whose value is $f$ on $D^c$, converges a.e. 
through
a certain parabola that depends on $\phi$ and $\phi'$.  
Our result includes the case $\phi(\lambda)=\log(1+\lambda^{\alpha/2})$.
Our proofs use both the probabilistic and analytic methods. In particular, the Poisson kernel estimates recently obtained in \cite{KK} are essential to our approach.  

\end{abstract}
\vspace{.6truein}



\section{Introduction}\label{sec-intro}

The classical Fatou theorem states that 
if $f \in L^p(\R^{d-1})$ for $p \in [1, \infty]$, 
then the Poisson extension $u_f$ of $f$ on the upper half-space     
has a nontangential limit a.e. on $\R^{d-1}$. 
It is also proved in \cite{L} that the nontangential approach is 
sharp.

Presently, non-local operators and their potential theory
have been extensively studied owing to their importance both in theories and applications. In particular, in \cite{BY, BY2}, the Fatou-type theorem for harmonic functions with respect to the operator
$\Delta^{\alpha/2}=-(-\Delta)^{\alpha/2}$ was discussed. 
Note that 
R. F. Bass and D. You 
\cite{BY} showed that the precise analogue of the Fatou theorem for harmonic functions with respect to 
$\Delta^{\alpha/2}$ is not true. Thus, in this case, it is necessary to state certain assumptions
 related to exterior functions to prove the existence of limits at the boundary.
Under certain $L^p$-H\"older continuity assumptions,    
it is shown in \cite{BY, BY2} that the Poisson extension 
with respect to $\Delta^{\alpha/2}$  has a nontangential limit a.e. on the upper half-space and Lipschiz domains, respectively.

Among many generalizations of the Fatou theorem,  
it has been proved that, under various types of assumptions on the boundary functions, the nontangential 
approach can be relaxed (see \cite{Do, NRS, NS, Su}). 
The boundedness of modified maximal operators 
has been an essential tool to prove this type of results. Recently, Y. Mizuta \cite{Mi2} 
applied analytic tools to the Poisson kernel of $\Delta^{\alpha/2}$ in the half space and showed 
that under the same assumption as that in \cite{BY}, 
a regular harmonic function (Poisson extension) with respect to $\Delta^{\alpha/2}$ in the half space  
has tangential limits a.e. if the approaching region is a certain parabola depending on $\alpha$.

The purpose of this paper is to investigate possible tangential approaching regions for a large class of the non-local operator 
$\phi(\Delta):=-\phi(-\Delta)$ in $C^{1,1}$ open sets. Our result extends the main result in  \cite{Mi2}. 
In \cite{Mi2}, the explicit Poisson kernel formula for $\Delta^{\alpha/2}$ in the upper half-space 
played a key role in proving the main result. 
However,
in general, 
it is not possible to derive an explicit Poisson kernel formula 
for $\phi(\Delta)$ in $C^{1,1}$ open sets. 
Fortunately, in a recent study \cite{KK}, we have obtained sharp two-sided estimates on the Poisson kernel for $\phi(\Delta)$ in bounded $C^{1,1}$ open sets under mild assumptions on $\phi$. 
In this paper, we use the upper bound of the Poisson kernel for $\phi(\Delta)$ near the boundary in \cite{KK} and 
show that the regular harmonic function with respect to $\phi(\Delta)$, 
which is the local $L^p$-H\"older continuous function of order 
$\bt$ on $D^c$ with $ p\in(1,\infty]$ and $\bt>1/p$, 
converges a.e. through a certain parabola. 
In our results, the approaching region depends on $\phi$ and $\phi'$. 
Nonetheless, our approaching region is always sufficiently wide to contain a Stolz open set.
See Remark \ref{r:region} to see how wide our approaching region is.

Before stating our procedure and the main result,
we introduce the following notations:
We use ``$:=$" to denote a definition, which
is  read as ``is defined to be". We denote $a \wedge b := \min \{ a, b\}$,
$a \vee b := \max \{ a, b\}$, 
 and $B(x,r):=\{y\in\Rd:|x-y|<r\}$.
For a set $W$ in $\R^d$, $\overline{W}$ and $W^c$ denote the closure and complement of $W$ in $\R^d$, respectively. For any open set $V$, we denote by $\delta_V (x)$,
the distance of a point $x$ to the boundary of $V$, i.e.,
$\delta_V(x)=\text{dist} (x,\partial V)$.
We often denote point $z=(z_1, \dots, z_d)\in \bR^d$ as
$(\wt z, z_d)$ with $\wt z \in \bR^{d-1}$.  
Since we consider tangential limits, we always assume that $d \ge 2$.

A smooth function $\phi:(0,\infty)\to [0,\infty)$ is called a Bernstein function if $(-1)^n\phi^{(n)}\le 0$ for every positive integer $n$. It is well-known that every Bernstein function with $\phi(0+)=0$ has the form
\begin{equation}\label{e:bernstein-function}
\phi(\lambda)=b \lambda +\int_{(0,\infty)}(1-e^{-\lambda t})\, \mu(\rd t)\, ,\quad \lambda >0\, ,
\end{equation}
where $b\ge 0$, and $\mu$ is a measure on $(0,\infty)$ satisfying $\int_{(0,\infty)}(1\wedge t)\, \mu(\rd t)<\infty$. 
$\mu$ is called the L\'evy measure of $\phi$.
(See \cite{SSV}.)

By concavity,
every Bernstein function $\phi$ satisfies
\begin{equation}\label{e:Berall}
\phi(t\lambda)\le
\lambda\phi(t)
\qquad \text{ for all }
\lambda \ge 1, t >0.
\end{equation}
Thus, 
$\lambda \mapsto \phi(\lambda)/\lambda$ 
is decreasing, and therefore, 
\begin{equation}\label{e:Berall1}
\lambda \phi'(\lambda)\le
\phi(\lambda)
\qquad \text{ for all }
\lambda >0.
\end{equation}
These simple properties of $\phi$ will be used several times in this paper.

In \cite{KM2}, the following 
conditions on the Bernstein function
$\phi$ are considered: 
Since we always assume that $d \ge 2$, here, we state the conditions only for $d \ge2$.

\noindent {\bf (A-1)}
$\phi$ is a complete Bernstein function, i.e., the L\' evy measure $\mu$ of $\phi$ has a 
completely monotone density $\chi$, i.e., $(-1)^{n}\chi^{(n)}\ge 0$ for every non-negative integer $n$.    
    
\noindent {\bf(A-2)} $\phi(0+)=0$ and $\lim_{\lambda \to \infty}\phi(\lambda)=\infty$.

\noindent
{\bf(A-3)}
        There exist constants $\sigma>0$, $\lambda_0>0$, and
$\delta \in (0, 1]$ such that
\begin{equation*}\label{eq:as-1}
  \frac{\phi'(\lambda t)}{\phi'(\lambda)}\leq
\sigma\,t^{-\delta}\ \text{ for all
}\ t\geq 1\ \text{ and }\ \lambda\geq\lambda_0\,.
\end{equation*}
 
 \noindent {\bf (A-4)}
If $d = 2$, we assume that there are $
\sigma_0>0$ and  
$\delta_0 \in (0,2\de)$ 
such that 
\begin{equation*}\label{e:new23}
  \frac{\phi'(\lambda t)}{\phi'(\lambda)}\geq
\sigma_0\,t^{-
\delta_0}\ \text{ for all
}\ t\geq 1\ \text{ and }\ \lambda\geq\lambda_0.
\end{equation*}
 
\noindent {\bf(A-5)}
  If the constant  $\delta$ in {\bf(A-3)} satisfies  $0<\delta \le
\frac{1}{2}$, then we assume that there exist constants $\sigma_1>0$ and
$\delta_1 \in [\delta, 1 )$ 
such that
\begin{equation*}\label{eq:as-2}
  \frac{\phi(\lambda t)}{\phi(\lambda)}\geq
\sigma_1\,t^{1-\delta_1}\ \text{ for all
}\ t\geq 1\ \text{ and }\ \lambda\geq\lambda_0\,.
\end{equation*}

\noindent {\bf(A-6)} There exist a $\te>0$ such that 
\begin{equation*}
\int_{0}^{\te}\frac{\ld^{d/2-1}}{\phi(\ld)}\rd \ld<\infty.
\end{equation*}
From {\bf(A-3)}, we get $b=0$ in \eqref{e:bernstein-function} by letting $t\to\infty$.
 From \cite[Lemma 2.2]{KM}, {\bf(A-3)} also implies that for every $\ep>0$, there exists $c=c(\ep)>0$
such that
\begin{align}\label{e:Hup}
\frac{\phi(\ld x)}{\phi(\ld)}\le cx^{1-\de+\ep}\qquad \text{for all}\;\; x\ge1\;\;\text{and}\;\;\ld\ge\ld_0.
\end{align}
See Example \ref{e:phi} for examples of  $\phi$ satisfying the assumptions {\bf(A-1)}--{\bf(A-6)}.

By Bochner's functional calculus, one can define the operator
$\phi(\Delta)$ on
$C_b^2(\R^d)$, which is the collection of bounded $C^2$ functions in $\RR^d$ with bounded derivatives.
Analytically, harmonic function $u$ for $\phi(\Delta)$ solves  
 $\phi(\Delta)u=0$ on $D$
in the distributional sense (see \cite{C1}). 
Since we use a probabilistic method, we formulate harmonic functions for $\phi(\Delta)$ using the L\'evy process 
corresponding to 
$\zt \mapsto\phi(|\zt|^2)$. 
Let 
$X=(X_t, \P_x)_{t\ge 0, x \in \R^d}$ be a rotationally
symmetric L\'evy process with a characteristic exponent
$\phi(|\zt|^2)$, 
that is, 
$$
\E_x\left[e^{i\zt\cdot(X_t-X_0)}\right]=e^{-t\phi(|\zt|^2)}
\quad \quad \mbox{ for every } x\in \R^d \mbox{ and } \zt\in \R^d.
$$
The infinitesimal generator of $X$ is $\phi(\Delta)$, i.e., 
$\phi(\Delta)u(x)=\lim_{t \to 0} t^{-1}( \E_x[u(X_t)]-u(x))$. 
For an open set $D$, let $\tau_D:=\inf\{t>0: \, X_t\notin D\}$. 
Now, we give the probabilistic definition of a (regular) harmonic function. 

\begin{definition}\label{d:har}
(1) A function 
$u:\bR^d\to \bR$ 
is said to be harmonic in an open set $D\subset \bR^d$ with respect to $X$ if
for every open set $B$ whose closure is a compact subset of $D$, $\bE_x[|u(X_{\tau_B})|]<\infty$ and 
$u(x)=\bE_x[u(X_{\tau_B})]$ for every $x\in B$.

\noindent
(2) A function 
$u:\bR^d\to \bR$ 
is said to be 
regular harmonic in an open set $D\subset \bR^d$ 
with respect to $X$ if $\bE_x[|u(X_{\tau_D})|]<\infty$ and 
$u(x)=\bE_x[u(X_{\tau_D})]$ for every $x\in D.$
\end{definition}
Clearly, a regular harmonic function in $D$ is harmonic in $D$ by the strong Markov property. 
Note that, by the Harnack inequality proved in \cite{KM}, 
under assumptions {\bf(A-1)}--{\bf(A-3)}  
the condition $\bE_x[|u(X_{\tau_B})|]<\infty$ for all $x \in D$ 
is equivalent to $\bE_{x_0}[|u(X_{\tau_B})|]<\infty$ for some $x_0 \in D$.

Now, we recall some function spaces related to our exterior functions.

\begin{definition}\label{cc}
Suppose $p \in (1, \infty]$. 
 \begin{description}
\item{(1)}
 $\Ld^{p}_{\bt}(\bR^d)$ is the space of $L^p$-H\"older continuous functions 
of order $\bt$ defined on $\bR^d$, 
i.e., 
$\bar f\in\Ld^{p}_{\bt}(\bR^d)$
means that $\bar f\in L^p(\bR^d)$ and there exists a constant $c>0$ such that \begin{equation}\label{e:holder}
\|\bar f(\cdot+y)-\bar f(\cdot)\|_{L^p(\bR^d)}\le c|y|^{\bt}\qquad\text{for all}\;\;y\in\bR^d.
\end{equation}
\item{(2)}
$\Ld_{\bt,loc}^{p}(\Dc)$ is the collection of functions $f$ such that $ f$ is defined on $\Dc$ 
and for each $\xi\in\partial D$,
there exists $\eta>0$ depending on $\xi$ such that $f$ agrees on $\Dc\cap B(\xi,\eta)$ with a function in 
$\Ld^{p}_{\bt}(\bR^d)$.
\end{description}
\end{definition}
Note that functions in $\Ld_{\bt,loc}^{p}(\Dc)$ may not be bounded (cf., \cite{St, BY, BY2}).

\begin{definition}\label{c11}
An open set $D$ 
in $\bR^d$ is said to be
$C^{1,1}$ if there exist  $R, \Lambda>0$ such that the following 
holds: for 
every $\xi\in\partial D$, there exist 
\begin{enumerate}
\item a $C^{1,1}$-function $\Gm=\Gm_{\xi}: \bR^{d-1}\to \bR$
satisfying $\Gm(0)=0$, $ \nabla\Gm (0)=(0, \dots, 0)$, $\| \nabla
\Gm \|_\infty \leq \Lambda$, $| \nabla \Gm (x)-\nabla \Gm (w)|
\leq \Lambda |x-w|$, for $x,w\in \bR^{d-1}$ and 
\item an orthonormal coordinate system $CS_{\xi}:y=( \wt y, \, y_d)$ with origin at $\xi$
such that 
$$B(\xi, R )\cap D= \{y=(\wt y, y_d) \in B(0, R) \mbox{
in } CS_{\xi}: y_d > \Gm ( \wt y) \}.$$ 
\end{enumerate}
The pair $( R, \Lambda)$ is called the $C^{1,1}$ characteristics of the open set $D$.
\end{definition}

For $\gm, a >0$, an $C^{1,1}$ open set $D$, 
and $\xi\in\partial D$,   
define
\begin{align}\label{s:appreg}
T_{\gm, \phi, a}(\xi)=T_{\gm, \phi, a, D}(\xi):=\left\{x\in D:|x-\xi|^{\gm+d}\phi(|x-\xi|^{-2})^{1/2}\le a
\frac{\de_D(x)^{d+2}\phi(\de_D(x)^{-2})^{3/2}}{\phi'(\de_D(x)^{-2})}\right\},
\end{align}
and $T_{\gm, \phi}(\xi):=T_{\gm, \phi, 1}(\xi)$.

Now, we state our theorem.
We use $\sH^s$ to denote the $s$-dimensional Hausdorff measure on $\Rd$ 
and for a measurable subset $W\subset \bR^d$, $|W|$ denotes the Lebesgue measure of $W$ in $\bR^d$. 
\begin{thm}\label{t:main}
 Suppose that
$p\in (1, \infty]$ and $\bt >1/p$. 
Let $X=(X_t, \P_x)_{t\ge 0, x \in \R^d}$ be a rotationally
symmetric L\'evy process with the characteristic exponent
$\phi(|\zt|^2)$ 
such that
the assumptions {\bf(A-1)}--{\bf(A-6)} hold and 
$\de$ in {\bf(A-3)} satisfies $1/p< \de \le 1$.
Suppose that $D$ is a $C^{1,1}$ open set with characteristic $(R, \Ld)$
and that 
$f\in\Ld_{\bt,loc}^{p}(\Dc)$ satisfies $\bE_{x_0}[ |f(X_{\tau_D})|]<\infty$ 
for some $x_0\in D$. 
Then, for $0<\gm<\beta-1/p$ and $a>0$, there exists a measurable subset 
$E\subset \partial D$ with  $\sH^{d-1}(E)=0$ such that
$u_{f}(x)=\bE_x[ f(X_{\tau_D})]$ 
has a finite limit along $T_{\gm, \phi, a}(\xi)$
for every $\xi\in\partial D\sms E$. 
Furthermore,
$$
\lim_{T_{\gm, \phi, a}(\xi)\ni x\to\xi}\left|u_f(x)-\lim_{r\to0+}\frac{1}{|B(\xi, r)\setminus D|}\int_{ B(\xi, r)\setminus D}f(y)\rd y\right|=0,\quad \text{for all}\;\;\xi\in\partial D\sms E.
$$
\end{thm}

Note that when $\phi(\lambda)=\lambda^{\alpha/2}$ and $D$ is the upper half-space $H:=\{x=(\wt x,x_d)\in\bR^d:x_d>0\}$,
our approaching region $T_{\gm, \phi, 2/\alpha}(\xi)$ is simply 
$\{x\in H:|x-\xi|^{1+\gm/(d-\alpha/2)}\le x_d\}$. 
Thus, our Theorem \ref{t:main} covers the result stated in 
\cite{Mi2}.
 
 Since the positive constant $a$ in \eqref{s:appreg} plays no special role in the proof, for convenience, we will only consider
$T_{\gm, \phi}(\xi)=T_{\gm, \phi, 1}(\xi)$.

\begin{example}\label{e:phi}\rm
Here are some examples of $\phi$ that satisfy the above assumptions {\bf(A-1)}--{\bf(A-6)}.
\begin{itemize}
\item $\phi(\ld)=\ld^{\alpha/2},\quad \alpha\in(0,2)$;
\item $\phi(\lambda)=(\lambda+\lambda^\alpha)^\kp,\quad \alpha, \kp\in (0, 1)$;
\item $\phi(\ld)=(\ld + m^{2/\alpha})^{\alpha/2}-m,\quad \alpha\in(0,2), \;m>0, \;d>2$;
\item $\phi(\ld)=\ld^{\af/2}+\ld^{\kp/2},\quad 0\le\kp<\alpha\in(0,2)$;
\item
$\phi(\lambda)=\lambda^{\alpha/2}(\log(1+\lambda))^{\kp},
\quad\alpha\in (0, 2)$, 
$\kp\in (-\alpha/2, 1-\alpha/2)$;
\item $\phi(\ld)=\log(1+\ld^{\alpha/2}),\quad \alpha\in(0,2], \;d>\alpha$;
\item $\phi(\ld)=\log(1+(\ld + m^{2/\alpha})^{\alpha/2}-m),\quad \alpha\in(0,2), \;m>0,\; d>2$.
\end{itemize}
\end{example}

\begin{remark}\label{r:region}
{\rm 
From \eqref{e:Berall1}, we see that $T_{\gm,\phi}(\xi)$ contains
\begin{align*}
T'_{\gm,\phi}(\xi):=\left\{x\in D:|x-\xi|^{\gm+d}\phi(|x-\xi|^{-2})^{1/2}\le
\de_D(x)^{d}\phi(\de_D(x)^{-2})^{1/2}\right\}.
\end{align*}
Moreover, $T'_{\gm, \phi}(\xi)$ contains the Stolz open set
\begin{align*}
S_M(\xi):=\{x\in D: |x-\xi|\le M\de_D(x), |x-\xi|<M^{-d/\gamma}\}
\end{align*}
for $M>1$. In fact, since $r^d\phi(r^{-2})^{1/2}$ is increasing by \eqref{e:Berall}, for $x\in S_M(\xi)$, we have 
$$|x-\xi|^d\phi(|x-\xi|^{-2})^{1/2}\le M^d\de_D(x)^d\phi(M^{-2}\de_D(x)^{-2})^{1/2}\le M^d\de_D(x)^d\phi(\de_D(x)^{-2})^{1/2}.$$
Thus, for $x\in S_M(\xi)$,
$$|x-\xi|^{\gm+d}\phi(|x-\xi|^{-2})^{1/2}\le 
|x-\xi|^{\gm}M^d\de_D(x)^d\phi(\de_D(x)^{-2})^{1/2} <
\de_D(x)^d\phi(\de_D(x)^{-2})^{1/2}.$$ 
We conclude that $S_M(\xi)\subset T'_{\gm,\phi}(\xi)\subset T_{\gm,\phi}(\xi)$.

On the other hand, our approaching region $T_{\gm,\phi}(\xi)$ can be strictly larger than $T'_{\gm,\phi}(\xi)$.
For example,  
when $D$ is the upper half-space $H$
and $\phi(\lambda)=\log(1+\lambda^{\alpha/2})$
where $\af\in(0,2]$, $d>\af$, 
we have 
\begin{align*}
T_{\gm, \phi}(\xi)
=&\left\{x\in H: |x-\xi|^{\gm+d}\{\log(1+|x-\xi|^{-\af})\}^{1/2}\le (2/\af)
(1+x_d^{\af})x_d^d\{\log(1+x_d^{-\af})\}^{3/2}\right\}\\
\supset& \left\{x\in H: |x-\xi|^{\gm+d}\{\log(1+|x-\xi|^{-\af})\}^{1/2}\le x_d^d\{\log(1+x_d^{-\af})\}^{3/2}\right\},
\end{align*}
while
\begin{align*}
T'_{\gm, \phi}(\xi)=\left\{x\in H:|x-\xi|^{\gm+d}\{\log(1+|x-\xi|^{-\af})\}^{1/2}\le x_d^d\{\log(1+x_d^{-\af})\}^{1/2}\right\}.
\end{align*}
}
\end{remark}

The remainder of this paper is organized as follows: In Section 2, we recall 
some basic facts on the Bernstein functions and corresponding L\'evy processes. Then, we recall 
the result in \cite{KK}, which is essential 
in  proving Theorem \ref{t:main}. Section 3 consists of key lemmas that hold for Lipschitz open sets. Using the results
in Section 2 and 3, we prove Theorem \ref{t:main} in Section 4. 
 
In this paper, we use the following convention: The values of
the constants  $R, \Ld, \ld_0, \de$ remain the same
throughout this paper, while $c, c_0, c_1, c_2, \ldots$  represent constants
whose values are unimportant and may change. All constants are positive finite numbers.
The constants $c_0, c_1, c_2, \ldots$ are labeled again in the statement and proof of
each result. The dependence of constant $c$ on dimension $d$ is not mentioned explicitly. 
For $x\in\bR^d$, $r>0$, and a set $W \subset \bR^d$, $B_W(x,r):=B(x,r)\cap W$.

\section{Preliminaries}

First, we recall certain essential facts about our L\'evy process $X$ that we will use later.
Then, since our proof considerably relies on the Poisson kernel estimates in \cite{KK},
we will also recall facts related to Poisson kernel and the result in \cite{KK}. 

Let $B=(B_t:\, t\ge 0)$ be a Brownian motion in
$\R^d$ whose infinitesimal generator is $\Delta$ (our Brownian motion $B$ runs
at twice the usual speed), and let $S=(S_t:\, t\ge 0)$ be a subordinator (non-negative increasing L\'evy process
in $\bR$ with $S_0=0$) 
independent of $B$ whose Laplace exponent is $\phi$, 
i.e., 
$$\E[\exp\{-\lambda S_t\}]=\exp\{-t\phi(\lambda)\}, \qquad \lambda >0.$$
It is well-known that the Laplace exponents of subordinators are always Bernstein functions. 
The L\'evy process $X=(X_t:\, t\ge 0)$ whose characteristic exponent is 
$\zt \mapsto \phi(|\zt|^2)$ 
can be defined 
by $X_t=B_{S_t}$ and it is also called a
subordinate Brownian motion. 
For example, a
rotation invariant 
$\alpha$-stable process is a subordinate Brownian motion with $\phi(\ld)=\ld^{\alpha/2}$.

For the remainder of this paper,  we will always assume that $\phi$ is a Bernstein function satisfying {\bf(A-1)}--{\bf(A-6)}.
Recall that  $\phi$  has the representation in \eqref{e:bernstein-function}. 
Since $b=0$ in \eqref{e:bernstein-function} by {\bf (A-3)},  $X$ is a pure jump process.

The L\'evy measure of $X$ has the density 
$x \mapsto j(|x|)$, 
where
$$
j(r)=\int_0^{\infty} (4\pi t)^{-d/2} e^{-r^2/(4t)}\, \mu(\rd t)\, 
$$
and $\mu$ is the L\'evy measure of $\phi$ (or of $S$). 
The infinitesimal generator of $X$ is $\phi(\Delta)$, which is an integro-differential operator of the type
$$
\phi(\Delta)u(x) =\int_{\R^d}\left(u(x+y)-u(x)-\nabla u(x) \cdot y {\mathbf 1}_{\{|y|\le 1\}}\right)\, j(|y|)\, \rd y .
$$
$X$ has a transition density $p(t,x,y)$ given by
\begin{align*}
p(t,x,y)=\int_{0}^{\infty}(4\pi s)^{-d/2}\exp\left(-\frac{|x-y|^2}{4s}\right)\bP(S_t\in \rd s).
\end{align*}
Recall that $X$ is said to be transient if $\bP_0(\lim_{t\to\infty}|X_t|=\infty)=1$. From the Chung-Fuchs type criterion of the transience of $X$, 
{\bf (A-6)} is equivalent to the transience of 
$X$ 
(see \cite[(2.9)] {KM}). 
Thus, we can define the Green function $G(x,y)$ by
\begin{align}
G(x,y)=g(|x-y|)=\int^{\infty}_{0}p(t,x,y)\rd t. \label{e:gg}
\end{align}
From \eqref{e:gg}, we see that $g$ is decreasing.
Under our assumptions {\bf(A-1)}--{\bf(A-6)},  $g(r)$ and $j(r)$ enjoy the following estimates
(see \cite{KM}): 
for every $M>0$, there exists $c=c(M)>0$ such that 
\begin{align}\label{H:2}
c^{-1}\frac{\phi'(r^{-2})}{r^{d+2}\phi(r^{-2})^2} \le g(r) \le c  
\frac{\phi'(r^{-2})}{r^{d+2}\phi(r^{-2})^2}, \qquad 0<r \le M,
\end{align}
and 
\begin{align}\label{H:3}
 \quad c^{-1}\frac{\phi'(r^{-2})}{r^{d+2}}\le  j(r)\le c\frac{\phi'(r^{-2})}{r^{d+2}}, \qquad 0<r \le M.
\end{align}

For any open subset $U$ in $\bR^d$, we use $G_U(x,y)$ to denote
the Green function of the process $X$ in $U$, which can be defined as $G_U(x,y)=G(x,y)-\bE_x[G(X_{\tau_U},y)]$.
For each fixed $z_0 \in U$, the function  $G_{U}(\cdot,z_0)$ is the non-negative regular harmonic function for $X$ in $U\sms \overline{B(z_0,\ep)}$ for every $\ep>0$ and it vanishes on $\R^d \setminus  U$.

Now, we define the Poisson kernel by 
\begin{equation*}\label{PK}
K_U(x,z)\,:=  \int_{U}
G_U(x,y) j(|y-z|) \rd y, \qquad (x,z) \in U \times
  {\overline U}^c.
\end{equation*}
Then, by the result of Ikeda and Watanabe (see \cite[Theorem 1]{IW}),
for
any open subset $U$ and  
every non-negative measurable function $f$,
\begin{align*}
\E_x\left[f(X_{\tau_U});\,X_{\tau_U-} \not= X_{\tau_U}  \right]
=\int_{\overline{U}^c} K_U(x,z)f(z)\rd z.
\end{align*}

\begin{definition}\label{Lips}
An open set $D$ in $\bR^d$ is said to be a Lipschitz
open set if there exist a localization radius $R_{\textrm{Lip}}>0$ and a constant
$\Lambda_{\textrm{Lip}} >0$ such that the following 
holds: for
every $\xi\in\partial D$, there 
exist 
\begin{enumerate}
\item a Lipschitz function $\psi=\psi_\xi:\bR^{d-1}\to \bR$ 
satisfying $\psi (0)= 0$, $| \psi (x)- \psi (y)|\le \Lambda_{\textrm{Lip}}|x-y|$, and 
\item an orthonormal coordinate system $CS_\xi:y=(\wt y, \, y_d)$ with its origin at
$\xi$ such that
$$
B(\xi, R_{\textrm{Lip}})\cap D=\{ y=(\wt y, y_d)\in B(0, R_{\textrm{Lip}}) \mbox{ in } CS_\xi: y_d
> \psi (\wt y) \}.
$$
\end{enumerate}
The pair $(R_{\textrm{Lip}}, \Lambda_{\textrm{Lip}})$ is called the characteristics of the
Lipschitz open set $D$.
\end{definition}

Since $X$ is a rotationally invariant pure jump L\' evy process, 
for every Lipschitz open set $D$,
$\bP_x(X_{\tau_D-} \not= X_{\tau_D})=1
$
(see \cite{M, Sz}). 
Thus, for every Lipschitz open set $D$ and every measurable function $f$ on $\Rd$, which satisfies
$
\int_{{\overline D}^c}K_D(x_0,z) | f(z)|\rd z<\infty$ for some $x_0 \in D$, 
$u_{ f}$ defined in Theorem \ref{t:main} has the following integral representation: 
\begin{equation}\label{e:rharmonic}
u_f(x)=\bE_x \left[ f(X_{\tau_D}):X_{\tau_D} \in  {\overline D}^c \right]
=\int_{{\overline D}^c}K_D(x,z)  f(z)\rd z, \qquad x\in D. 
\end{equation}
Clearly, any regular harmonic function $u$ in a Lipschitz open set $D$, whose value on $D^c$ is $f$, is written as $u_f$.

Furthermore, when $U$ is a bounded $C^{1,1}$ open set, 
we see from \cite[Theorem 1.3]{KK} that 
\begin{align}
&c^{-1}\frac{\phi(\delta_U(z)^{-2})^{1/2}}{\phi(\delta_U(x)^{-2})^{1/2}\phi(|x-z|^{-2})(1+\phi(\delta_U(z)^{-2})^{-1/2})}j(|x-z|)
\nonumber\\
&\le  K_U(x,z) \le
c\,\frac{\phi(\delta_U(z)^{-2})^{1/2}}{\phi(\delta_U(x)^{-2})^{1/2}\phi(|x-z|^{-2})(1+\phi(\delta_U(z)^{-2})^{-1/2})}j(|x-z|), \qquad (x,z) \in U \times
  {\overline U}^c.\label{estK}
\end{align}
We will use the upper bound in \eqref{estK} for $|x-z| <2$.

\section{Analysis on Lipschitz open set}
Recall that we assume {\bf(A-1)}--{\bf(A-6)}. 
In this section, we prove some results that hold on Lipschitz
open sets. We will use these results in Section 4.

Throughout this section, we fix the Lipschitz open set $D$ with characteristics $(R_{\textrm{Lip}}, \Lambda_{\textrm{Lip}})$.
Without loss of generality, we assume that
$R_{\textrm{Lip}}<1$.
Note that $D$ can be
unbounded and disconnected. For every
$\xi\in \partial D$ and $ x \in B(\xi, R_{\textrm{Lip}})$, 
we define the vertical distance
$$
\rho_\xi (x) := x_d-\psi_\xi (\wt x)\, ,
$$
where $(\wt x, x_d)$ are the coordinates of $x$ in $CS_\xi$.
Then,
\begin{align}\label{e:compara}
\de_D(x)\le|\rho_\xi(x)|\le(1+\Lambda_{\textrm{Lip}})\de_D(x), \qquad \xi\in \partial D,\;  x \in B(\xi, R_{\textrm{Lip}}).
\end{align}

Recall that $\ld_0$ and $\delta$ are the constants in {\bf(A-3)}. 
\begin{lemma}\label{l:intphi}
For all $q\in[1,1/(1-\de))$, and $M\ge 1$, 
there exists a constant $c=c(q, \de, \Lambda_{\textrm{Lip}}, M)>0$  
such that 
for every $\xi\in \partial D$,  
$s \le R_{\textrm{Lip}}/2$,
and $r \le (2M)^{-1} (R_{\textrm{Lip}} \wedge \ld_0^{-1/2})$,
\begin{align}\label{e:intphi}
\int_{\{(\wt y, y_d) \text{ in } CS_{\xi}: |\wt y|<s, |\rho_\xi(y)|< M r\}}\phi(\de_D(y)^{-2})^{q/2}\rd y 
\le cr s^{d-1}{\phi(r^{-2})^{q/2}}.
\end{align}
\end{lemma}
\pf 
First, since $\phi$ is increasing, by \eqref{e:compara}, the left-hand side of\eqref{e:intphi} is less than or equal to 
\begin{align}\label{e:pdiwq}
\int_{\{(\wt y, y_d) \text{ in } CS_{\xi}: |\wt y|<s, |\rho_\xi(y)|< M r\}} \phi\Big((1+\Lambda_{\textrm{Lip}})^2|\psi_{\xi}(\wt y)-y_d|^{-2}\Big)^{q/2} \rd y .
\end{align}
Using the assumption 
$q<1/(1-\de)$, 
choose $\ep=\ep(\de,q) \in (0, \delta+1/q-1)$.
By the change of variable $t=\rho_\xi(y)/M$, the fact that $\phi$ is increasing, 
\eqref{e:Berall}, and  \eqref{e:Hup},
\begin{align*}
&\int_{\{(\wt y, y_d) \text{ in } CS_{\xi}: |\wt y|<s, |\rho_\xi(y)|< M r\}} \phi\Big((1+\Lambda_{\textrm{Lip}})^2|\psi_{\xi}(\wt y)-y_d|^{-2}\Big)^{q/2} \rd y \nn\\
&=
\phi((1+\Lambda_{\textrm{Lip}})^2 M^{-2} r^{-2})^{q/2}\int_{\{(\wt y, y_d) \text{ in } CS_{\xi}: |\wt y|<s, |\rho_\xi(y)|< M r\}} \left(\frac{\phi ((1+\Lambda_{\textrm{Lip}})^2|\psi_{\xi}(\wt y)-y_d |^{-2})}
{\phi((1+\Lambda_{\textrm{Lip}})^2 M^{-2} r^{-2})}
\right)^{q/2}\rd y\nn\\
&\le\;2(1 \vee (1+\Lambda_{\textrm{Lip}})^q M^{-q})\phi(r^{-2})^{q/2}
\int_{\{|\wt y|<s\}}\int_{0}^{r}\left(\frac{\phi((1+\Lambda_{\textrm{Lip}})^2 M^{-2}t^{-2})}{\phi((1+\Lambda_{\textrm{Lip}})^2 M^{-2}r^{-2})}\right)^{q/2}M\rd t\rd \wt y \nn\\
&\le\;c_1\phi(r^{-2})^{q/2} s^{d-1}\int_{0}^{r}\left(\frac{r}{t}\right)^{(1-\de+\ep)q}\rd t, \end{align*}
which is less than or equal to $c_2 \phi(r^{-2})^{q/2} s^{d-1} r$ since $(1-\de+\ep)q<1$.
Combining this and \eqref{e:pdiwq},  we have proved the lemma.\qed

A direct consequence of Lemma \ref{l:intphi} is the following: 
\begin{cor}\label{c:intphi}
For all $\delta\in(0,1]$, 
there exists a constant $c=c(\de, \Lambda_{\textrm{Lip}})>0$ such that for every   
$r \le (2+2\Lambda_{\textrm{Lip}})^{-1}(R_{\textrm{Lip}} \wedge \ld_0^{-1/2})$ 
and all $\xi\in \partial D$, 
$$
\int_{B(\xi,r)}\phi(\de_D(y)^{-2})^{1/2}\rd y \le cr^{d}\phi(r^{-2})^{1/2}.
$$
\end{cor}

Recall that for $x\in\bR^d$, $r>0$, and a set $W \subset \bR^d$, $B_W(x,r)=B(x,r)\cap W$.

\begin{lemma}\label{l:base}
Suppose that 
$f\in \Ld_{\bt,loc}^{p}(\Dc)$
and $0<\gm<\beta-1/p$.
 If $\de$ in {\bf(A-3)} satisfies $\de>1/p$, 
then $\sH^{d-1}(E(\gm))=0$, where 
\begin{align*}
E(\gm):=\left\{\xi\in \partial D:\limsup_{r\to0+}
\int_{B_{\Dc}(\xi,r)}\int_{B_{\Dc}(\xi,r)}\frac{\phi(\de_D(y)^{-2})^{1/2}}{r^{2d+\gm}\phi(r^{-2})^{1/2}}| f(y)-f(z)|\rd y \rd z >0\right\}.
\end{align*}
\end{lemma}

\pf 
 Using the cardinality, we can choose $\xi_i\in \partial D$ and $\eta_i=\eta_i(\xi_i)>0$ such that
there exists $\bar f_i\in \Ld^{p}_{\bt}(\bR^d)$ satisfying $f=\bar f_i$ on $\Dc\cap B(\xi_i,\eta_i)$, and that, 
$\partial D\subset \cup_{i\in\bN}B(\xi_i, a_i)$, where $a_i=(\RL\wedge \eta_i)/4$.
Let 
$E_i(\gm)= E(\gm)\cap B(\xi_i,a_i)$ and $n_i\in \bN$ be the largest number such that 
$n_i\le \log_2 (a_i^{-1} \vee \sqrt{\ld_0})+\log_2(2+2\Lambda_{\textrm{Lip}})+1$.
For $M>0$ and $n\ge n_i$, 
set
\begin{align*}
&E_i(\gm,M,n)\nn\\
&=\left\{\xi\in B_{\partial D}(\xi_i,a_i):\int_{B_{\Dc}(\xi,2^{-n})}\int_{B_{\Dc}(\xi,2^{-n})}
\frac{\phi(\de_D(y)^{-2})^{1/2}}{2^{-n(2d+\gm)}\phi(2^{2n})^{1/2}}| \bar f_i(y)- \bar f_i(z)|\rd y \rd z>\frac1M\right\}.
\end{align*}
Since
$$
E(\gm) \subset \bigcup_{i=0}^\infty E_i(\gm)= \bigcup_{i=0}^\infty \bigcup_{M=1}^\infty \left(\bigcap^{\infty}_{k=n_i}\bigcup^{\infty}_{n=k}E_i(\gm, M,n)\right),$$
it suffices to show that 
\begin{align}\label{BC}
\sum^{\infty}_{n= n_i}\sH^{d-1} (E_i(\gm,M,n))<\infty,
\end{align} 
which implies $\sH^{d-1}(\bigcap^{\infty}_{k=n_i}\bigcup^{\infty}_{n=k}E_i(\gm, M,n))=0$ by Borel-Cantelli Lemma.
Throughout the remainder of the proof, we fix $M$ and $i$ and assume that $n\ge n_i$.

Let $\psi=\psi_{\xi_i}$ and $CS=CS_{\xi_i}$ be the Lipschitz function and the orthonormal coordinate system in Definition \ref{Lips}. We will use this coordinate system $CS$ below so that $\xi_i=0$.
For $\xi:=(\wt\xi, \psi(\wt\xi))\in  B_{\partial D}(0,a_i)$ in $CS$, define
\begin{align}
h_n(\wt\xi):=\int_{B_{\Dc}(\xi, 2^{-n})}\int_{B_{\Dc}(\xi, 2^{-n})}
\phi(\de_D(y)^{-2})^{1/2}|\bar f_i(y)-\bar f_i(z)|\rd y\rd z.\nn
\end{align} 
Then, by using the area formula (see, for example, \cite[Section 3.3.4]{EG}),
\begin{align}\label{e:measure}
\sH^{d-1}(E_i(\gm,M,n))
=&\int_{B_{\partial D} (0,a_i)}{\bf 1}_{E_i(\gm,M,n)}(\wt\xi, \psi(\wt\xi))\rd\sH^{d-1}(\wt\xi, \psi(\wt\xi))\nn\\
\le&\;M2^{n(2d+\gm)}\phi(2^{2n})^{-1/2}\int_{B_{\partial D} (0,a_i)}h_n(\wt\xi)\rd\sH^{d-1}(\wt\xi, \psi(\wt\xi))\nn\\
\le&\;M2^{n(2d+\gm)}\phi(2^{2n})^{-1/2}\int_{|\wt\xi|<a_i }h_n(\wt\xi)(1+|\nabla\psi(\wt\xi)|^2)^{1/2}\rd\wt\xi\nn\\
\le&\;(1+\LL^2)^{1/2}M2^{n(2d+\gm)}\phi(2^{2n})^{-1/2}\int_{|\wt\xi|<a_i}h_n(\wt\xi)\rd\wt\xi.
\end{align}

When $p\in(1,\infty)$,
by H\"older's inequality,
\begin{align}\label{e:hpart}
&\int_{|\wt\xi|<a_i}h_n(\wt\xi)\rd\wt\xi \nonumber \\
\le &\; \int_{|\wt\xi|<a_i}\int_{B(0,2^{-n})}\int_{B(0,2^{-n})}
\phi(\de_D(\xi+y)^{-2})^{1/2}
|\bar f_i(\xi+y)-\bar f_i(\xi+z)|\rd y \rd z\rd \wt\xi\nn\\
\le &\; \left(|B(0,2^{-n})|\int_{|\wt\xi|<a_i}\int_{B(0,2^{-n})}
\phi(\de_D(\xi+y)^{-2})^{q/2}\rd y \rd \wt\xi\right)^{1/q}\nn\\
&\qquad\times\left(\int_{|\wt\xi|<a_i}\int_{B(0,2^{-n})}\int_{B(0,2^{-n})}
|\bar f_i(\xi+y)-\bar f_i(\xi+z)|^{p}\rd y \rd z\rd \wt\xi\right)^{1/p}=: \;I\times II,
\end{align}
where $1/q:=1-1/p$.

By Fubini's theorem,
\begin{align}\label{e:Iq}
I^{q}\le c_1 2^{-nd}\int_{|\wt y|<2^{-n}}\int_{|y_d|<2^{-n}}\int_{|\wt\xi|<a_i}
\phi(\de_D(\xi+y)^{-2})^{q/2}\rd \wt\xi\rd y_d\rd\wt y,
\end{align}
while using $|a+b|^p\le 2^{p-1}(|a|^{p}+|b|^{p})$, the symmetry, and Fubini's theorem,
\begin{align}\label{e:II2}
II^{p}&=\int_{|\wt\xi|<a_i}\int_{B(0,2^{-n})}\int_{B(0,2^{-n})}
|\bar f_i(\xi+y)-\bar f_i(\xi+z)|^{p}\rd y \rd z\rd \wt\xi\nn\\
&\le 2^{p-1}\int_{|\wt\xi|<a_i}\int_{B(0,2^{-n})}\int_{B(0,2^{-n})}
|\bar f_i(\xi+y)-\bar f_i(\xi+y+z)|^{p}\rd y \rd z\rd \wt\xi\nn\\
&\quad+2^{p-1}\int_{|\wt\xi|<a_i}\int_{B(0,2^{-n})}\int_{B(0,2^{-n})}
|\bar f_i(\xi+y+z)-\bar f_i(\xi+z)|^{p}\rd y \rd z\rd \wt\xi\nn\\
&\le 2^p \int_{B(0,2^{-n})}\int_{|\wt y|<2^{-n}}\int_{|y_d|<2^{-n}}\int_{|\wt\xi|<a_i}
|\bar f_i(\xi+y+z)-\bar f_i(\xi+y)|^{p}\rd \wt\xi\rd y_d \rd \wt y\rd z.
\end{align}

Let $w=(\wt w, w_d):=(\wt\xi + \wt y,\psi(\wt\xi)+y_d)=\xi+y$. If $|\wt y|< 2^{-n}$ and
 $|\psi(\wt w-\wt y)-w_d|=|y_d|<2^{-n}$, then
$$|w_d-\psi(\wt w)|\le |w_d-\psi(\wt w-\wt y)|+|\psi(\wt w-\wt y)-\psi(\wt w)|
\le 2^{-n}+\LL|\wt y|\le(1+\LL)2^{-n}.$$
Therefore, for $|\wt y|< 2^{-n}$, 
\begin{align}  
&\{w\in \bR^d:|\wt w-\wt y|<a_i,|\psi(\wt w-\wt y)-w_d|<2^{-n}\}\nn\\
&\subset \{w\in \bR^d:|\wt w|<2a_i, |\psi(\wt w)-w_d|<(1+\LL)2^{-n}\}=:{Q_n}.\label{subset}
\end{align}
Using this and Lemma \ref{l:intphi}, the inner two integrals in \eqref{e:Iq} are bounded as 
\begin{align}\label{e:I3}
\int_{|y_d|<2^{-n}}\int_{|\wt\xi|<a_i}\phi(\de_D(\xi+y)^{-2})^{q/2}\rd \wt\xi \rd y_d
\le\int_{{Q_n}}\phi(\de_D(w)^{-2})^{q/2}\rd w_d \rd \wt w 
\le c_2 2^{-n} \phi(2^{2n})^{q/2}.
\end{align}
Further, from \eqref{e:holder}, 
the inner two integrals in \eqref{e:II2} are bounded as 
\begin{align}\label{e:II3}
\int_{|y_d|<2^{-n}}\int_{|\wt\xi|<a_i}|\bar f_i(\xi+y+z)-\bar f_i(\xi+y)|^{p}\rd \wt\xi\rd y_d
\le \int_{Q_n}|\bar f_i(w+z)-\bar f_i(w)|^{p}\rd w
\le c_3^p 2^{-n\beta p}.
\end{align}

Thus,  \eqref{e:Iq}, \eqref{e:I3}  imply 
$I\le c_4 2^{-2nd/q} \phi(2^{2n})^{1/2}$
and 
\eqref{e:II2}, \eqref{e:II3} imply 
$II\le  c_5 2^{-2nd/p}2^{-n( \beta -1/p)}.$
From this and \eqref{e:hpart}, we obtain
\begin{align}\label{e:hxi}
\int_{|\wt\xi|<a_i}h_n(\wt\xi)\rd\wt\xi\le \;c_6  2^{-n(2d+ \beta -1/p)}\phi(2^{2n})^{1/2}.
\end{align}
Now, we conclude from \eqref{e:measure} and \eqref{e:hxi} that
$\sH^{d-1}(E_i(\gm,M,n))
\le\;c_7 2^{-n(\beta-1/p-\gm)},$ 
which implies \eqref{BC} 
since $\beta-1/p-\gm>0$.

When $p=\infty$, simply by \eqref{e:holder} and Corollary  \ref{c:intphi},
\begin{align*}
h_n(\wt\xi) \le c_{8}(2^{-n+1})^\bt\int_{B_{\Dc}(\xi, 2^{-n})}\int_{B_{\Dc}(\xi, 2^{-n})}
\phi(\de_D(y)^{-2})^{1/2}\rd y\rd z 
\le c_{9}2^{-n(2d+\bt)} \phi(2^{2n})^{1/2}.
\end{align*} 
Therefore, by \eqref{e:measure}, 
$\sH^{d-1}(E_i(\gm,M,n))
\le c_{10} 2^{-n(\bt-\gm)},$ 
which yields \eqref{BC} 
since $\bt-\gm>0$. 

\qed

\begin{lemma}\label{l:base2}
Suppose that 
$f\in \Ld_{\bt,loc}^{p}(\Dc)$
and $0<\gm<\beta-1/p$.
Let
\begin{align*}
F(\gm)=\left\{\xi\in \partial D:\limsup_{r\to0+}r^{-2d-\gm}
\int_{B_{\Dc}(\xi,r)}\int_{B_{\Dc}(\xi,r)}|f(y)-f(z)|\rd y \rd z>0\right\},
\end{align*}
then $\sH^{d-1}(F(\gm))=0$.
\end{lemma}
\pf  The proof of this lemma is the same as that of Lemma \ref{l:base}. 
In fact, using the same 
$a_i$, $\bar f_i$, $n_i$, and coordinate system 
in the proof of Lemma \ref{l:base}, for $n \ge n_i$, we define
\begin{align*}
F_i(\gm,M,n)
:=\left\{\xi\in  B_{\partial D}(0,a_i):2^{n(2d+\gm)}
\int_{B_{\Dc}(\xi,2^{-n})}\int_{B_{\Dc}(\xi,2^{-n})}|\bar f_i(y)- \bar f_i(z)|\rd y \rd z>\frac1M\right\}.
\end{align*}

When $p\in(1,\infty)$, 
by H\"older's inequality we have  
\begin{align*}
\sH^{d-1}(F_i(\gm,M,n))
\le&\;(1+\LL^2)^{1/2}M2^{n(2d+\gm)}\int_{|\wt\xi|<a_i}\int_{B_{\Dc}(\xi,2^{-n})}\int_{B_{\Dc}(\xi,2^{-n})}
|\bar f_i(y)- \bar f_i(z)|\rd y \rd z\rd\wt\xi\\
\le& \;(1+\LL^2)^{1/2}M2^{n(2d+\gm)} \left(\int_{|\wt\xi|<a_i}\int_{B(0,2^{-n})}\int_{B(0,2^{-n})}
\rd y \rd z\rd \wt\xi\right)^{1/q}\nn\\
&\quad\times\left(\int_{|\wt\xi|<a_i}\int_{B(0,2^{-n})}\int_{B(0,2^{-n})}
|\bar f_i(\xi+y)-\bar f_i(\xi+z)|^{p}\rd y \rd z\rd \wt\xi\right)^{1/p}.
\end{align*}
By following the proof of Lemma \ref{l:base} line by line, 
we see that
\begin{align*}
\left(\int_{|\wt\xi|<a_i}\int_{B(0,2^{-n})}\int_{B(0,2^{-n})}
\rd y \rd z\rd \wt\xi\right)^{1/q}
= c_{1} \left(a_i^{d-1}2^{-2nd}\right)^{1/q}, 
\end{align*}
and
\begin{align*}
\left(\int_{|\wt\xi|<a_i}\int_{B(0,2^{-n})}\int_{B(0,2^{-n})}
|\bar f_i(\xi+y)-\bar f_i(\xi+z)|^{p}\rd y \rd z\rd \wt\xi\right)^{1/p}
\le c_{2} (2^{-nd}2^{-n(d-1)}2^{-n\beta p})^{1/p}.
\end{align*}
Therefore, 
\begin{align*}
\sH^{d-1}(F_1(\gm,M,n))
\le\; c_{3}2^{n(2d+\gm)}2^{-2nd}2^{n/p}2^{-n\bt}=\;c_{3}2^{-n(\beta-1/p-\gm)}.
\end{align*}
The assertion for $p\in(1,\infty)$ follows from this. 

The proof for $p=\infty$ is also similar. Therefore, we skip the proof. 
\qed

For a locally integrable function $h$ on $\bR^d$ and bounded measurable set $U\subset \bR^d$, 
we define its integral mean over the region $U$ by
$\fint_{U}h(y)\rd y=\frac{1}{|U|}\int_{U}h(y)\rd y.$

The proof of the following lemma is taken from \cite{Mi2}. However, for the reader's convenience, we state the details of the proof.
\begin{lemma}\label{l:existA}
Let $f\in \Ld_{\bt,loc}^{p}(\Dc)$, and $\de$ in {\bf(A-3)} satisfies $\de>1/p$. 
Then,
$$A(\xi):=\lim_{r\to0+}\fint_{ B_{\Dc}(\xi, r)}f(y)\rd y$$
exists and is finite for $\sH^{d-1}$-a.e. 
$\xi\in\partial D$.  
Moreover, for $0<\gm<\beta-1/p$,
\begin{align}
\lim_{r\to0+}r^{-d-\gm}\phi(r^{-2})^{-1/2}\int_{B_{\Dc}(\xi,r)}\phi(\de_D(y)^{-2})^{1/2}|f(y)-A(\xi)|\rd y=0 \label{e:existA}
\end{align}
for $\sH^{d-1}$-a.e. 
$\xi\in\partial D$. 
\end{lemma}
\pf
For simplicity, define $ A(\xi,r)=\fint_{B_{\Dc}(\xi,r)}f(y)\rd y$.
Then, for $r\le t\le 2r$,
\begin{align*}
|A(\xi,t)-A(\xi,r)|\le c_1r^{-2d}\int_{B_{\Dc}(\xi,2r)}\int_{B_{\Dc}(\xi,2r)}|f(y)-f(z)|\rd y\rd z.
\end{align*}
Lemma \ref{l:base2} gives
$\lim_{r\to0+}r^{-\gm}|A(\xi,2r)-A(\xi,r)|=0$
for $\sH^{d-1}$-a.e. 
$\xi\in\partial D$. 
This implies that $A_{\infty}(\xi):=\lim_{n\to\infty}A(\xi,2^{-n})$ exists and 
$\lim_{k\to\infty}2^{k\gm}\left\{A(\xi,2^{-k+1})-A_{\infty}(\xi)\right\}=0.$ 
Therefore, for $r\le 2^{-k+1}\le 2r$,
\begin{align*}
r^{-\gm}|A_{\infty}(\xi)-A(\xi,r)|
&\le r^{-\gm}|A_{\infty}(\xi)-A(\xi,2^{-k+1})|+r^{-\gm}|A(\xi,2^{-k+1})-A(\xi,r)|\\
&\le 2^{k\gm}|A_{\infty}(\xi)-A(\xi,2^{-k+1})|+r^{-\gm}|A(\xi,2^{-k+1})-A(\xi,r)| 
\longrightarrow 0 \quad \text{as} \quad r\to 0. 
\end{align*}
Thus, $A(\xi)$ exists and is finite $\sH^{d-1}$-a.e. 
$\xi\in\partial D$. 
Further, for $\sH^{d-1}$-a.e. 
$\xi\in\partial D$, 
\begin{align}\label{e:lim}
\lim_{r\to0+}r^{-\gm}|A(\xi)-A(\xi,r)|=0.
\end{align}
By Corollary \ref{c:intphi},  
\begin{align*}
&\;r^{-d-\gm}\phi(r^{-2})^{-1/2}\int_{B_{\Dc}(\xi,r)}\phi(\de_D(y)^{-2})^{1/2}|f(y)-A(\xi)|\rd y \\
\le&\;r^{-d-\gm}\phi(r^{-2})^{-1/2} \int_{B_{\Dc}(\xi,r)}\phi(\de_D(y)^{-2})^{1/2}|f(y)-A(\xi,r)|\rd y
+c_2r^{-\gm}|A(\xi)-A(\xi, r)|\\
\le&\;c_3 r^{-2d-\gm}\phi(r^{-2})^{-1/2}\int_{B_{\Dc}(\xi,r)}\int_{B_{\Dc}(\xi,r)}\phi(\de_D(y)^{-2})^{1/2}
|f(y)-f(z)|\rd z\rd y
+c_2r^{-\gm}|A(\xi)-A(\xi, r)|, 
\end{align*}
which tends to zero as $r\to 0$ for $\sH^{d-1}$-a.e.
$\xi\in\partial D$ 
by \eqref{e:lim} and Lemma \ref{l:base}. 
Hence, we have proved  \eqref{e:existA}.
\qed

\section{Proof of Theorem \ref{t:main}}

For any $C^{1,1}$ open set $D$ with characteristic $(R,\Ld)$, it is well-known that
(see, e.g., \cite[Lemma 2.2]{Song})
there exists $L=L(R,\Ld,d)>0$ such that for every $\xi\in \partial D$ and $r\le (R \wedge 1)$, one can
obtain a $C^{1,1}$ open set $U(\xi,r)$ with characteristic $(r(R \wedge 1)/L, \Ld L/r)$ such that
\begin{align}\label{e:bddc1}
D\cap B(\xi,r/2)\subset U(\xi,r)\subset D\cap B(\xi,r).
\end{align}

First, we record a lemma, which is
a consequence of the main results of \cite{KM2}. Although simple, it is important in this paper.

\begin{lemma}\label{l:bL2}
Let $D$ be a $C^{1,1}$ open set with characteristic $(R,\Ld)$. 
Suppose that $\xi_0 \in \partial D$, $r_0>0$, and 
$u:\bR^d\to \bR$ 
is a non-negative regular harmonic function in $D\cap B(\xi_0,r_0)$  with respect to $X$ vanishing on $D^c\cap B(\xi_0,r_0)$. 
Then, $u$ vanishes continuously on $(\partial D)\cap B(\xi_0,r_0/2)$.  
\end{lemma}

\pf Fix $\xi \in (\partial D)\cap B(\xi_0,r_0/2)$ and  
let $0<r < 1 \wedge R \wedge (r_0/2)$. We will show that $u$ vanishes continuously on $(\partial D)\cap B(\xi ,r /8)$, which clearly implies the lemma. 

Choose a bounded $C^{1,1}$ open set $U=U(\xi, r)$ as in \eqref{e:bddc1}.
Let $z_0$ be a point in $U \sms B(\xi,r/4)$. 
Then, by \cite[Proposition 2.4 and Theorem 7.1]{KM2}, the $G_{U}(\cdot ,z_0)$ vanishes 
continuously on $\partial U \supset (\partial D)\cap B(\xi,r/4) $. Moreover, $ x \mapsto G_{U}(x ,z_0)$ is a regular harmonic function in $D\cap B(\xi,r/4)$  with respect to $X$. 
Thus, by 
the boundary Harnack principle  
{\cite[Theorem 5.6(i)]{KM2}}, for a fixed $x_0\in B(\xi,r/8)\cap D$ and $x \in  D\cap B(\xi,r/8)$, 
\begin{align*}
u(x)
& \le c u(x_0){G_U(x_0, z_0)^{-1} }G_{U}(x, z_0)
\longrightarrow 0 \quad \text{as}\;\;x\to (\partial D)\cap B(\xi,r/8).
\end{align*}
This completes the proof.
\qed

Before we prove the main result, we observe an inequality. 
Recall that $g(r)$ defined in \eqref{e:gg} is decreasing. Using this fact and the estimates in  \eqref{H:2}, we obtain that there exists $c>0$ such that
\begin{align}
\frac{\phi'(t^{-2})}{\phi(t^{-2})t^{d+2}}\le c \frac{\phi'(s^{-2})}{\phi(s^{-2})s^{d+2}}, \quad   s\le t\le 2.   \label{e:Gdec}
\end{align}

Now, we prove our theorem.

\vspace{0.5cm}

\noindent{\bf  Proof of Theorem \ref{t:main}.}
Without loss of generality, we assume that $R<1$.
By the cardinality, we can choose $\xi_i\in \partial D$ and $\eta_i=\eta_i(\xi_i)>0$ such that
there exists $f_i\in \Ld^{p}_{\bt}(\bR^d)$ satisfying $f=f_i$ on 
$\Dc\cap B(\xi_i,\eta_i)$, and that,  
$\partial D\subset \cup_{i\in\bN}B(\xi_i,\eta_i/8)$.
Without loss of generality, we let $\eta_i\le R$.
Since $\partial D$ is a countable union of $B_{\partial D}(\xi_i, \eta_i/8)$,
it suffices to show that for $\sH^{d-1}$-a.e. 
$\xi\in P_1:=B_{\partial D}(\xi_1, \eta_1/8)$, 
$u_{f}$ has a limit along $T_{\gm, \phi}(\xi)$.

Choose a $C^{1,1}$ open set $U=U(\xi_1,R)$ with characteristic $(R^2/L, \Ld L/R)$ as \eqref{e:bddc1} so that 
$P_1=B_{\partial U}(\xi_1, \eta_1/8)$. 
Note that 
\begin{align*}
u_{f}(x)=\bE_{x}[f(X_{\tau_D})]
&= \bE_{x}[f(X_{\tau_D});\tau_{U}<\tau_D]+\bE_{x}[f(X_{\tau_{U}});\tau_{U}=\tau_D]\\
&= \bE_{x}[f(X_{\tau_D});\tau_{U}<\tau_D]+\bE_{x}[f(X_{\tau_{U}})]
-\bE_{x}[f(X_{\tau_{U}});\tau_{U}<\tau_D].
\end{align*}
By using the strong Markov property at $\tau_{U}$, 
$\bE_{x}[f(X_{\tau_D});\tau_{U}<\tau_D]$ 
and $\bE_{x}[f(X_{\tau_{U}});\tau_{U}<\tau_D]$
are non-negative regular harmonic functions  in $U$ (and thus, in $B(\xi_1,R/2)\cap D$)
with respect to $X$ and vanish on $D^c\cap B(\xi_1,R/2)$. 
Thus, by Lemma \ref{l:bL2}, the limits of $u_{f}(x)$ 
and $\bE_{x}[f(X_{\tau_{U}})]$ (if exist) are the same when $x$ goes to a point $\xi\in P_1$. 
Therefore, it suffices to show that the limit 
$\lim_{T_{\gm, \phi}(\xi) \ni x \to \xi} \bE_x[f(X_{\tau_{U}})]$ 
exists for
$\sH^{d-1}$-a.e $\xi\in P_1$.

By Lemma \ref{l:existA}, for $\sH^{d-1}$-a.e. 
$\xi\in P_1$, 
we have that 
$A(\xi)=\lim_{r\to0+}\fint_{ B_{\Uc}(\xi, r)}f(y)\rd y$
exists and is finite and that \begin{align}
\lim_{r\to 0+}r^{-\gm-d}\phi(r^{-2})^{-1/2}\int_{B_{\Uc}(\xi,r)}\phi(\delta_U(y)^{-2})^{1/2}|f(y)-A(\xi)|\rd y =0 \label{e:final1}
\end{align}
holds.
For the remainder of the proof, we 
fix a  $\xi\in  P_1$ 
and show that \begin{align}
\lim_{T_{\gm, \phi}(\xi)\ni x\to\xi}|\bE_x[f(X_{\tau_{U}})]-A(\xi)|=0. \label{e:final}
\end{align}

Let $\ep>0$ be given.
By \eqref{e:final1},
there exists $r_0<(1\wedge (R/4))/2$ such that for every $0<r<2r_0$,
\begin{align}\label{e:ball}
\int_{B_{\Uc}(\xi,r)}\phi(\delta_U(y)^{-2})^{1/2}|f(y)-A(\xi)|\rd y <\ep r^{\gm+d}\phi(r^{-2})^{1/2}.
\end{align}
Note that 
$B(\xi,2r_0)\subset B(\xi_1,R/2)$.  
Let
\begin{align*}
u_1(x)=\bE_x[|f(X_{\tau_{U}})|;
X_{\tau_{U}}\in \bR^d\sms\{\overline{U}\cup B(\xi,r_0)\}],\;\;
u_2(x)=\bP_x(X_{\tau_{U}}\in \bR^d\sms\{\overline{U}\cup B(\xi,r_0)\}).
\end{align*}
Then, $u_1, u_2$ are non-negative regular harmonic functions in $U\cap B(\xi,r_0)$
with respect to $X$, and they vanish on 
$U^c\cap B(\xi,r_0)$. 

On the other hand, since $U$ is a bounded $C^{1,1}$ open set, 
we have the following Poisson kernel estimates by \eqref{estK} and \eqref{H:3}:
\begin{align}\label{e:Kupper2}
K_{U}(x,y) \le
c_1\,\frac{\phi(\delta_U(y)^{-2})^{1/2}}{\phi(\delta_U(x)^{-2})^{1/2}\phi(|x-y|^{-2})}
\frac{\phi'(|x-y|^{-2})}{|x-y|^{d+2}}
,\quad \text{for}\;\;x\in U, 
\;y\in B(\xi,r_0)\setminus{\overline U}.
\end{align}
Thus, by \eqref{e:rharmonic} and \eqref{e:Kupper2}, we have 
\begin{align}
|\bE_x[f(X_{\tau_{U}})]-A(\xi)|&\le 
c_1  \int_{B_{\overline{U}^c}(\xi, r_0)} \frac{\phi(\delta_U(y)^{-2})^{1/2}}{\phi(\delta_U(x)^{-2})^{1/2}\phi(|x-y|^{-2})}
\frac{\phi'(|x-y|^{-2})}{|x-y|^{d+2}}|f(y)-A(\xi)|\rd y\nn\\
&\quad+u_1(x)+|A(\xi)|u_2(x). \label{e:newpp} 
\end{align}

Since $|x-y| \ge \delta_U(x)$ for $y \in \bR^d\sms \overline{U}$,
by \eqref{e:Gdec} and \eqref{e:ball}, 
for $x \in B(\xi,r_0/8) \cap T_{\gm, \phi}(\xi)$,
\begin{align}\label{e:nbdry}
&
\int_{
B_{\overline{U}^c}(\xi, 2|x-\xi| ) } 
\frac{\phi(\delta_U(y)^{-2})^{1/2}}{\phi(\delta_U(x)^{-2})^{1/2}}
\left(\frac{\phi'(|x-y|^{-2})}{\phi(|x-y|^{-2})|x-y|^{d+2}}\right)|f(y)-A(\xi)|\rd y\nn\\
&\le c_2\frac{\phi'(\delta_U(x)^{-2})}{\delta_U(x)^{d+2}\phi(\delta_U(x)^{-2})^{3/2}}
\int_{
B_{\overline{U}^c}(\xi, 2|x-\xi| ) } \phi(\delta_U(y)^{-2})^{1/2}|f(y)-A(\xi)|\rd y\nn\\
&\le c_2 2^{d+\gm} \frac{\phi'(\delta_U(x)^{-2}) |x-\xi|^{d+\gm}\phi(2^{-2}  |x-\xi|^{-2})^{1/2}}{\delta_U(x)^{d+2}\phi(\delta_U(x)^{-2})^{3/2}}\ep \le c_2 2^{d+\gm}\ep.
\end{align}

When $2|x-\xi|\le  |\xi -y|$, we have 
$|x-y| \ge |y-\xi|-|\xi-x| \ge |\xi-y|/2$. Thus, by \eqref{e:Gdec} and \eqref{e:Berall1},
on $\{y \in  \bR^d\sms \overline{U}: 2|x-\xi|\le  |\xi -y| <r_0 \}$,
\begin{align*}
&\frac{\phi'(|x-y|^{-2})}{|x-y|^{d+2}\phi(|x-y|^{-2})}
\le c_2 \frac{2^{d+2}\phi'(4|\xi-y|^{-2})}{|\xi-y|^{d+2}\phi(4|\xi-y|^{-2})}
\le   \frac{c_22^{d}}{|\xi-y|^{d}}\nn\\
&\le   \frac{c_22^{2d}}{2^d-1} \left(\frac{1}{|\xi-y|^{d}}   -\frac{1}{(2r_0)^{d}}\right)=\frac{c_2 d 2^{2d}}{2^d-1}
\int^{2r_0}_{{2 |x-\xi|}} {\bf 1}_{\{|y-\xi|<t\}}\frac{\rd t}{t^{d+1}}.
\end{align*}
Therefore, by the Fubini's theorem and \eqref{e:ball}, we have that for $x \in B(\xi,r_0/8)$,
\begin{align}
&\int_{\{y \in  \bR^d\sms \overline{U}: 2|x-\xi|\le  |\xi -y| <r_0 \}} 
\frac{\phi(\delta_U(y)^{-2})^{1/2}}{\phi(\delta_U(x)^{-2})^{1/2}\phi(|x-y|^{-2})}
\frac{\phi'(|x-y|^{-2})}{|x-y|^{d+2}}|f(y)-A(\xi)|\rd y\nn\\
& \le \frac{c_2 d 2^{2d}}{2^d-1} \phi(\delta_U(x)^{-2})^{-1/2}
\int^{2r_0}_{{2 |x-\xi|}}\left(\int_{B_{\overline U^c}(\xi,t)}
\phi(\delta_U(y)^{-2})^{1/2}|f(y)-A(\xi)|\rd y\right)\frac{\rd t}{t^{d+1}} \nn\\
&\le \frac{c_2 d 2^{2d}}{2^d-1}\ep 
\phi(\delta_U(x)^{-2})^{-1/2}\phi({2^{-2} |x-\xi|}^{-2})^{1/2}\int^{2r_0}_{0}t^{\gm-1}\rd t\nn\\
&\le \frac{c_2 d 2^{2d} }{(2^d-1)\gm}\ep  (2r_0)^{\gm}\le \frac{c_2 d 2^{2d} }{(2^d-1)\gm}\ep.\label{e:intermed}
\end{align}

Applying \eqref{e:nbdry} and  \eqref{e:intermed} to \eqref{e:newpp}, together with Lemma \ref{l:bL2} gives 
\begin{align*}
\limsup_{T_{\gm, \phi}(\xi)\ni x\to\xi}|\bE_x[f(X_{\tau_{U}})]-A(\xi)|
\le c_3\ep,
\end{align*}
where the constant $c_3>0$ is independent of $\ep$.  Since $\ep>0$ is arbitrary, 
we have proved the claim \eqref{e:final}.
\qed

\medskip

{\bf Acknowledgements.} 
We thank the referee for his (her) helpful comments on the first version of
this paper.

\vspace{.1in}
\begin{singlespace}

\end{singlespace}

{\bf Panki Kim}

Department of Mathematical Sciences and Research Institute of Mathematics,
Seoul National University,
Building 27, 1 Gwanak-ro, Gwanak-gu,
Seoul 151-747, Republic of Korea

E-mail: \texttt{pkim@snu.ac.kr}

\bigskip

{\bf Jaehoon Kang}

Department of Mathematical Sciences,
Seoul National University,
Building 27, 1 Gwanak-ro, Gwanak-gu,
Seoul 151-747, Republic of Korea

E-mail: \texttt{jaehnkang@gmail.com}

\bigskip

\end{document}